\documentclass[11pt]{article}
\usepackage{amsfonts,epsf,amsmath,amssymb}
\usepackage{tikz, caption, subcaption}
\usepackage[numbers,sort&compress]{natbib}

\newtheorem{theorem}{\bf Theorem}[section]
\newtheorem{corollary}[theorem]{\bf Corollary}

\newtheorem{lemma}[theorem]{\bf Lemma}
\newtheorem{proposition}[theorem]{\bf Proposition}
\newtheorem{conjecture}[theorem]{\bf Conjecture}

\tikzstyle{vertex}=[circle, draw, inner sep=0pt, minimum size=6pt]
\newcommand{\pch}{\chi_{\rho}}

\newcommand{\proof}{\noindent{\bf Proof.\ }}
\newcommand{\qed}{\hfill $\square$ \bigskip}
\newcommand{\smallqed}{{\tiny ($\Box$)}}

\begin{document}

\title{\bf Packing chromatic number versus chromatic and clique number}

\author{
\phantom{xxxx} \and Bo\v{s}tjan Bre\v{s}ar $^{a,b}$  \and Sandi Klav\v zar $^{a,b,c}$ \and \phantom{xx}
\and Douglas F. Rall $^{d}$ \and Kirsti Wash $^e$ }

\date{}

\maketitle

\begin{center}
$^a$ Faculty of Natural Sciences and Mathematics, University of Maribor, Slovenia\\
\medskip

$^b$ Institute of Mathematics, Physics and Mechanics, Ljubljana, Slovenia\\
\medskip

$^c$ Faculty of Mathematics and Physics, University of Ljubljana, Slovenia\\
\medskip

$^d$ Department of Mathematics, Furman University, Greenville, SC, USA\\
\medskip

$^e$ Department of Mathematics, Western New England University, Springfield, MA, USA\\
\end{center}

\begin{abstract}
The packing chromatic number $\pch(G)$ of a graph $G$ is the smallest integer $k$ such that the vertex set of $G$ can be partitioned into sets $V_i$,  $i\in [k]$, where each $V_i$ is an $i$-packing. In this paper, we investigate for a given triple $(a,b,c)$ of positive integers whether there exists a graph $G$ such that $\omega(G) = a$, $\chi(G) = b$, and $\chi_{\rho}(G) = c$. If so, we say that $(a, b, c)$ is realizable. It is proved that $b=c\ge 3$ implies $a=b$, and that triples $(2,k,k+1)$ and $(2,k,k+2)$ are not realizable as soon as $k\ge 4$. Some of the obtained results are deduced from the bounds proved on the packing chromatic number of the Mycielskian. Moreover, a formula for the independence number of the Mycielskian is given. A lower bound on $\pch(G)$ in terms of $\Delta(G)$ and $\alpha(G)$ is also proved.
\end{abstract}

\noindent {\bf Key words:} packing chromatic number; chromatic number; clique number; independence number; Mycielskian

\medskip\noindent
{\bf AMS Subj.\ Class:} 05C70, 05C15, 05C12


\section{Introduction}
A fundamental problem in graph coloring is the relation between the chromatic number $\chi(G)$ of a graph $G$ and its clique number $\omega(G)$. The construction of Mycielski provided examples of graphs that are triangle-free and have arbitrarily large chromatic number~\cite{myc-55}. Hence graphs with arbitrary clique number $k$ and chromatic number of an arbitrary size greater than $k$ could be constructed. In this paper, we ask similar questions involving the packing chromatic number by studying the existence of graphs $G$ with given $\omega(G)$, $\chi(G)$ and $\pch(G)$.

Given a graph $G$ and a positive integer $i$, an {\em $i$-packing} in $G$ is a subset $W$ of the vertex set of $G$ such that the distance between any two distinct vertices from $W$ is greater than $i$. This generalizes the notion of an independent set, which is equivalent to a $1$-packing. The {\em packing chromatic number} of $G$ is the smallest integer $k$ such that the vertex set of $G$ can be partitioned into sets $V_1,\ldots, V_k$, where $V_i$ is an $i$-packing for each $i\in [k] = \{1,\ldots, k\}$. This invariant is well defined on any graph $G$ and is denoted $\pch(G)$. More generally, for a nondecreasing sequence
$S=(s_1,\ldots,s_k)$ of positive integers, the mapping $c:V(G)\longrightarrow [k]$ is an {\em $S$-packing coloring} if for any $i$ in $[k]$ the set $c^{-1}(i)$ is an $s_i$-packing~\cite{goddard-2012}.

In particular, if $S=(1,\ldots,k)$, then $c:V(G)\longrightarrow [k]$ is called a {\em $k$-packing coloring}, which is the main concept in this paper.
The packing chromatic number was introduced in~\cite{goddard-2008} under the name broadcast chromatic number, and subsequently studied under the current name, see~\cite{argiroffo-2014, balogh-2017+, barnaby-2017, bkr-2007, bkr-2016, bkrw-2017a, bkrw-2017b, ekstein-2014, fiala-2010, fiala-2009, finbow-2010, jacobs-2013, korze-2014, lbe-2016, shao-2015, togni-2014, torres-2015}.

Clearly, in any graph $G$, $\omega(G)\le \chi(G)\le\pch(G)$, and the main question we are interested in is for which triples $(a,b,c)$, where $2\le a\le b\le c$, there exists a graph $G$ such that $\omega(G)=a$, $\chi(G)=b$ and $\pch(G)=c$. In this paper, we use the name {\em realizable triple} for a triple $(a,b,c)$ whose values are realized by some graph. As it turns out, the Mycielski construction is useful also in this study.
Recall that the {\em Mycielskian} $M(G)$ of a graph $G$ is the graph with the vertex set $V(G) \cup V'\cup \{w\}$, where $V' = \{x':\ x\in V(G)\}$, and the edge set $E(G) \cup \{xy':\ xy\in E(G)\}\cup \{wx':\ x'\in V'\}$. Well-known properties of this construction are that $\chi(M(G))=\chi(G)+1$ and $\omega(M(G))=\omega(G)$.

In studying realizable triples the following result from the seminal paper will be used several times.

\begin{proposition} {\rm (\cite[Proposition 2.1]{goddard-2008})}
\label{prp:indep}
If $G$ is a graph with order $n(G)$, then $\pch(G)\le n(G) - \alpha(G) + 1$, with equality if ${\rm diam}(G) = 2$.
\end{proposition}

In view of Proposition~\ref{prp:indep} and the usefulness of the Mycielskian in chromatic graph theory, in Section~\ref{sec:mycielskian-indep} we investigate the packing chromatic number and the independence number of the Mycielskian.  We present a formula for establishing $\alpha(M(G))$ in an arbitrary graph $G$; to the best of our knowledge, this has not yet been established in full generality, cf.~\cite{lin-2006}. The obtained formula is then applied to obtain various bounds on $\pch(M(G))$. We also show that the packing chromatic number of the Mycielskian $M(G)$ is at least two more than that of $G$.

In Section~\ref{sec:chi_vs_rho-packing} we first prove our main result asserting that $\chi(G)=\pch(G)$ implies that $\omega(G) = \chi(G)$. (It was proven in~\cite{goddard-2008} that $\chi(G)=\pch(G)$ implies that $\omega(G)\ge \chi(G)-2$.) In other words, $(a,b,b)$ is realizable only if $a=b$. Next, we prove that if $(a,b,c)$ is realizable, then $(a,b,d)$ is also realizable for any $d$ that is greater than $c$. If $k\ge 4$, we show that the triple $(2,k,k+2)$ is not realizable. On the other hand, by applying the Mycielskian operation several times, we infer that triples $(n,n+k,2^{k-1}(n+1)+1)$ are realizable for any $n\ge 2$ and any $k\ge 1$.

In the final section we present the following lower bound on the packing chromatic number of an arbitrary graph:
$$\pch(G)\ge \Delta(G) - \alpha(G) + 2,$$
which in a nice way complements the bound from Proposition~\ref{prp:indep} (here, $\Delta(G)$ stands for the maximum degree of vertices in $G$). Some of the graphs that attain this bound are used in presenting families of graphs, which realize $(a,a,a)$ for $a\ge 2$.

\section{Independence and packing chromatic number of the Mycielskian}
\label{sec:mycielskian-indep}

In view of Proposition~\ref{prp:indep} it is important to know the independence number of a graph while studying its packing chromatic number. In this section we also consider the independence number of the Mycielskian. Although the Mycielskian has been investigated by now from many points of view \cite{bara-2008,dos-2005,glg-2012,huch-1999,llz-2010}, it seems that for the independence number only sporadic results were obtained.

Setting ${\cal I}(G)$ to denote the set of independent sets of a graph $G$ (including the empty set), the independence number of the Mycielskian can be described as follows.

\begin{theorem}
\label{thm_mycielski-indep}
If $G$ is a connected graph, then
$$\alpha(M(G)) = \max_{S\in {\cal I}(G)} \{ 2|S| + |V(G)\setminus N[S]|\}\,.$$
\end{theorem}

\proof
Let $V(G) = \{v_1, \ldots, v_n\}$, so that the vertex set of $M(G)$ is $V(G) \cup V'\cup \{w\}$, where $V'=\{v_1', \ldots, v_n'\}$. Let $M = \max_{S\in {\cal I}(G)} \{ 2|S| + |V(G)\setminus N[S]|\}$.

Let $S\in {\cal I}(G)$. Set $S' = \{x':\ x\in S\}$ and $X' = \{y':\ y\in V(G)\setminus N[S]\}$. Since $S\in {\cal I}(G)$ we also have that $S\in {\cal I}(M(G))$. Since $S'\cup X'\subseteq V'$, we clearly have $S'\cup X'\in {\cal I}(M(G))$. It is also clear (since $S$ is independent) that there are no edges between $S$ and $S'$. Finally, since $x\in S$ has no neighbor in $V(G)\setminus N[S]$, there are also no edges between $S$ and $X'$. It follows that $S\cup S'\cup X'\in {\cal I}(M(G))$. Consequently, $\alpha(M(G))\ge M$.

To prove the reverse inequality, let $S$ be an arbitrary independent set of $M(G)$ with $|S| = \alpha(M(G))$. Let $S_G = S \cap V(G)$ and note that it is possible that $S_G = \emptyset$.
By the definition of $S$, the vertices of $V(G)\setminus S_G$ do not lie in $S$. Moreover, if $x\in V(G)\setminus S_G$ is adjacent to a vertex from $S_G$, then also $x'\in V'$ is not in $S$. But all the other vertices from $V'$ can lie in $S$ and since $S$ is a largest independent set, all these vertices do lie in $S$. Setting $Y' = \{y':\ y\in V(G)\setminus N[S]\}$ we thus have that $|S| = 2|S_G| + |Y'| = 2|S_G| + |V(G)\setminus N[S_G]|$. We conclude that $\alpha(M(G))\le M$.
\qed

The following consequence to Theorem~\ref{thm_mycielski-indep} was first proven in~\cite{lin-2006}. More precisely, it can be deduced from Theorems~4.1 and~4.2 of~\cite{lin-2006} by specializing to the case $m=1$ and by replacing the vertex cover number with the independence number.

\begin{corollary} {\rm \cite{lin-2006}}
\label{cor:upper-lower-mycielski}
If $G$ is a connected graph, then
$$2\alpha(G) \le \alpha(M(G)) \le n(G) + \alpha(G) - 1\,.$$
\end{corollary}

\proof
Let $S\in {\cal I}(G)$ and let $|S| = k$. Then $|V(G)\setminus N[S]|\le n(G) - k - 1$ and hence $2|S| + |V(G)\setminus N[S]| \le 2k + (n(G) - k - 1) = n(G) + k -1$. Since $k\le \alpha(G)$ it follows that $2|S| + |V(G)\setminus N[S]| \le n(G) + \alpha(G) -1$. The upper bound now follows from Theorem~\ref{thm_mycielski-indep}.

For the lower bound select $S\in {\cal I}(G)$ with $|S| = \alpha(G)$.
\qed

The upper bound of Corollary~\ref{cor:upper-lower-mycielski} can be improved as follows.

\begin{proposition}
\label{propo:upper-down-by-1}
If $G$ is a connected graph which is neither a complete graph nor a star, then $\alpha(M(G)) \le n(G) + \alpha(G) - 2$.
\end{proposition}

\proof
As usual, let $V(M(G)) = V(G) \cup V'\cup \{w\}$. Let $S\in {\cal I}(M(G))$ with $|S| = \alpha(M(G))$, and let $S' = S \cap V'$. If $w\in S$, then $S' = \emptyset$ and consequently $|S| \le 1 + \alpha(G) \le n(G) + \alpha(G) -2$, where the last inequality holds because $G$ is not complete and thus $n(G)\ge 3$. Hence we may assume in the rest that $w\notin S$.

If $|S'| = n(G)$, then necessarily $S' = S$ and hence the conclusion of the proposition holds because $G$ is not complete and thus $\alpha(G)\ge 2$.

Suppose next that $|S'| \le n(G) - 2$. Since $|S\cap V(G)|\le \alpha (G)$, we immediately get that $\alpha(M(G)) = |S| \le n(G) -2 + \alpha(G)$.

In the last case to be considered assume that $|S'| = n(G) - 1$. If $|S\cap V(G)|\le \alpha (G) - 1$, then the conclusion is clear (since $w\notin S$). Hence suppose that $|S\cap V(G)| = \alpha (G)$. Let $x'$ be the unique vertex of $V'$ that is not in $S'$. Since $G$ is connected and $x$ has a neighbor in $S'$, it
follows that $x\notin S$. Moreover, a vertex $y\in V(G)\setminus S$, $y\ne x$, would imply that $S$ is not independent. It follows that $|S\cap V(G)| = n(G) - 1$, that is, $\alpha(G) = |V(G)| - 1$. But this means that $G$ is a star.
\qed

Next, we turn our attention to the packing colorings of the Mycielskian.

\begin{theorem}
\label{thm:mycielski}
If $G$ is a connected graph with $n(G)\ge 2$, then $\pch(M(G))\ge \pch(G) + 2$, with equality if $G$ is complete.
\end{theorem}

\proof
Let $V(G) = \{v_1, \ldots, v_n\}$, so that the vertex set of $M(G)$ is $V(G) \cup V'\cup \{w\}$, where $V'=\{v_1', \ldots, v_n'\}$.

Consider first the case $G=K_n$. Since ${\rm diam}(M(K_n)) = 2$ and $\alpha(M(K_n)) = n$, Proposition~\ref{prp:indep} implies that $\pch(M(K_n)) = n(M(K_n)) - \alpha(M(K_n)) + 1 = n + 2$. Hence the result holds (with equality) for complete graphs. We may assume in the rest of the proof that $G$ is not complete, in particular, $n(G)\ge 3$.

Let $\pch(M(G)) = k$ and let $c$ be a $k$-packing coloring of $M(G)$. Note that $M(G)$ contains an induced $C_5$, hence $k\ge 4$. We distinguish the following cases.

\medskip\noindent
{\bf Case 1}: $c(w) = 1$. \\
In this case $c(v_i')\ne 1$ for $i\in [n]$. Moreover, $|c(V')| = |V'| = n(G)$. It follows that $\pch(M(G))\ge 1 + |V(G)|$. Since $G$ is not complete we have $\pch(G) \le |V(G)| - 1$. Consequently $\pch(M(G))\ge 1 + |V(G)| \ge \pch(G) + 2$.

\medskip\noindent
{\bf Case 2}: $c(w) = k$. \\
Let $\widetilde{c}$ be  a coloring defined on $V(G)$ as follows:
$$\widetilde{c}(v_i) =
   \left\{ \begin{array}{ll}
         c(v_i'); & c(v_i) = k-1\,,\\
         c(v_i); & \mbox{otherwise}.
    \end{array} \right.
$$
Note first that $\widetilde{c}: V(G) \rightarrow [k-2]$. Indeed, since $c(w)=k$ and ${\rm ecc}(w) = 2$, the color $k$ is not used by $\widetilde{c}$. In addition, since $d_{M(G)}(v_i,v_i') = 2$ we have that $c(v_i') \le k-2$ for any vertex $v_i$ with $c(v_i) = k-1$. We next claim that $\widetilde{c}$ is a packing coloring. Since the restriction of $c$ to $V(G)$ is a packing coloring, it suffices to show that setting $\widetilde{c}(v_i) = c(v_i') = \ell$ if $c(v_i) = k-1$, preserves the property of being packing coloring.

Assume that $\widetilde{c}(v_j) = \ell$ holds for some $j\ne i$. This holds because either (i) $c(v_j) = \ell$ or (ii) $c(v_j) = k-1$ and $c(v_j') = \ell$. Suppose first that (i) happened. Let $P$ be a shortest $v_i,v_j$-path in $G$, and let $x$ be the neighbor of $v_i$ on $P$. (It is possible that $x=v_j$.) Since $v_i'$ is in $M(G)$ adjacent to $x$, we have $d_{M(G)}(v_i', v_j) \le d_G(v_i,v_j)$. Since $c(v_i') = c(v_j) = \ell$, we have  $d_{M(G)}(v_i', v_j) > \ell$, which implies that $d_G(v_i,v_j) > \ell$ as required. Suppose next that (ii) holds, that is, $c(v_j) = k-1$ and $c(v_j') = \ell$. Since we also have $c(v_i') = \ell$ and as $d_{M(G)}(v_i', v_j') = 2$, we must have $\ell = 1$. But since $c(v_i) = c(v_j) = k-1$, we clearly have $d_G(v_i,v_j) > 1$.  We conclude that $\widetilde{c}$ is a packing coloring. Hence $\pch(G) \le k-2 = \pch(M(G)) - 2$.

\medskip\noindent
{\bf Case 3}: $2\le c(w) \le k-1$. \\
First let $\widetilde{c}$ be  a coloring defined on $V(G)$ as follows:
$$\widetilde{c}(v_i) =
   \left\{ \begin{array}{ll}
         c(v_i'); & c(v_i) = k\,,\\
         c(v_i); & \mbox{otherwise}.
    \end{array} \right.
$$
Since for $i\in [n]$ we have $d_{M(G)}(v_i,v_i') = 2$, we get that $\widetilde{c}: V(G) \rightarrow [k-1]$. Because $2\le c(w) \le k-1$ and ${\rm ecc}(w) = 2$, the color $c(w)$ is used only on vertex $w$. If $c(w) = k-1$, then $\widetilde{c}: V(G) \rightarrow [k-2]$ and as in Case 2 we see that $\widetilde{c}$ is a packing coloring.  Otherwise $\widetilde{c}: V(G) \rightarrow [k-1]\setminus \{c(w)\}$. Now let $\widehat{c}$ be the coloring of $G$ obtained from $\widetilde{c}$ by recoloring each vertex of color $k-1$ with color $c(w)$. Then $\widehat{c}: V(G) \rightarrow [k-2]$ is a required packing coloring.
\qed

The stars $K_{1,n}$, $n\ge 2$, form another family for which the equality is achieved in Theorem~\ref{thm:mycielski}. Indeed, it is easy to verify that $\pch(K_{1,n}) = 2$ and $\pch(M(K_{1,n})) = 4$. Another example is provided by the path $P_4$ for which we  have $\pch(P_4) = 3$ and $\pch(M(P_4)) = 5$.

On the other hand, the difference $\pch(M(G))-\pch(G)$ can be arbitrarily large. For example,  consider $K_{t,t}$, $t\ge 2$. The graph $M(K_{t,t})$ has diameter $2$, hence having in mind Proposition~\ref{prp:indep} and Theorem~\ref{thm_mycielski-indep},
\begin{eqnarray*}
\pch(M(K_{t,t}))-\pch(K_{t,t}) & = & n(M(K_{t,t}))-\alpha(M(K_{t,t}))+1-(t+1) \\
& = & 2n(K_{t,t})+1-2t+1-t-1 \\
& =& t+1\,.
\end{eqnarray*}

But we can bound $\pch(M(G))$ from the above as follows.

\begin{proposition}
\label{prp:mycielski-from-above}
If $G$ is a connected graph with $n(G)\ge 2$, then
$$\pch(M(G))\le \min\{ n(G) + 2, 2(n(G)-\alpha(G)+1)\}\,.$$
\end{proposition}

\proof
Again let the vertex set of $M(G)$ be $V(G) \cup V'\cup \{w\}$, where $V'=\{v':\ v\in V(G)\}$. Then $V'$ is an independent set of $M(G)$. Coloring vertices from $V'$ with color $1$ and every other vertex with a unique color greater than $1$ is a packing coloring using $n(G) + 2$ colors. Similarly, if $X$ is an independent set of $G$ with $|X|=\alpha(G)$, then $X\cup \{x':\ x\in X\}$ is an independent set of $M(G)$ of order $2\alpha(G)$, cf.\ Theorem~\ref{thm_mycielski-indep}. Proceeding as in the first case we find a packing coloring using $1 + 2(n(G)-\alpha(G)) + 1$ colors.
\qed

The lower bound of Theorem~\ref{thm:mycielski} coincides with the upper bound of Proposition~\ref{prp:indep} on complete graphs and on stars.
We note that if ${\rm diam}(G) = 2$, then ${\rm diam}(M(G)) = 2$ as well. Hence Proposition~\ref{prp:indep} implies:

\begin{corollary}
\label{cor:mycielski-diam-2}
If ${\rm diam}(G) = 2$, then $\pch(M(G)) = 2n(G) - \alpha(M(G)) + 2$.
\end{corollary}

Consider the following example. Let $G_{k,\ell}$, $k,\ell\ge 3$, be the graph obtained from the complete graph $K_k$ by selecting a vertex $x$ of $K_k$ and attaching $\ell$ pendant vertices to $x$. The diameter of  $G_{k,\ell}$ is 2. Using Theorem~\ref{thm_mycielski-indep} or directly we see that
$\alpha(M(G_{k,\ell})) = 2\ell +k - 1$. Hence Corollary~\ref{cor:mycielski-diam-2} implies that $\pch(G_{k,\ell}) = k + 3$.

We point out that the fact $\alpha(M(G_{k,\ell})) = 2\ell +k - 1$ demonstrates that there are graphs $G$ for which $\alpha(M(G))$ is arbitrarily far away from the lower bound given in Corollary~\ref{cor:upper-lower-mycielski}. 

\section{Realizing graphs of given clique, chromatic, and packing chromatic numbers}
\label{sec:chi_vs_rho-packing}
Given a sequence $(a, b, c)$ of positive integers, we say that {\it $(a, b, c)$ is realizable} if there exists a graph $G$ such that $\omega(G) = a$, $\chi(G) = b$, and $\chi_{\rho}(G) = c$, in which case we say $G$ {\em realizes $(a, b, c)$} and that $G$ is {\em an  $(a,b,c)$ graph}. We know that for any graph $G$, $\omega(G) \le \chi(G) \le \chi_{\rho}(G)$. Thus, $(a, b, c)$ must be a nondecreasing sequence in order for $(a,b,c)$ to be realizable. Furthermore, $\omega(G) \ge 2$ for any nontrivial graph so we only consider sequences where $a\ge 2$. For example, the only triangle-free $2$-chromatic graphs with packing chromatic number $2$ are stars. Therefore, the only graphs that realize $(2,2,2)$ are stars. A natural question to ask is whether a realizable sequence $(a,b,c)$ implies $(a,b,d)$ is realizable for any $d > c$. The following result answers this question in the affirmative.

\begin{lemma}\label{lem:abd}
If $(a, b, c)$ is realizable, then $(a, b, d)$ is realizable for every $d$, where $d>c$.
\end{lemma}

\proof
Let $G$ be a graph that realizes $(a, b, c)$. We first show that there exists a graph $G'$ which realizes $(a, b, d)$ for some $d>c$ and contains $G$ as a subgraph. Construct $G'$ by appending $c$ leaves to each vertex of $G$. Note that $G'$ has the same clique size and chromatic number as $G$. We claim that $\chi_{\rho}(G')=r>c$. To see this, let $f:V(G') \to [r]$ be a packing coloring of $G'$. Since $G$ is a subgraph of $G'$, we know that the restriction of $f$ to $V(G)$ is a packing coloring of $G$ so $r\ge c$. Moreover, if no vertex of $V(G)$ receives color $1$, then some vertex of $V(G)$ is assigned a color larger than $c$, for otherwise (by decreasing each color used on $V(G)$ by 1) it follows that $\pch(G)<c$, which is a contradiction. Hence, if no vertex of $V(G)$ receives color $1$, we have $r>c$. On the other hand, if there exists a vertex $v$ of $G$ such that $f(v)=1$, then the leaves appended to $v$ receive pairwise different colors. Thus, some leaf of $G'$ is given a color greater than $c$. It follows that $r>c$ and $(a,b,r)$ is realizable for some $r>c$.

Finally, to see that $(a, b, d)$ is realizable for all $d$, where $d>c$, we only need to show that $(a, b, c+1)$ is realizable. Indeed, pick a vertex $v$ of $G$ and append a leaf $w$ to $v$. Either $\pch(G+w) = \pch(G)$ or $\pch(G+w) = \pch(G) + 1$. If $\pch(G+w) = \pch(G)$, then continue appending leaves to $w$ until either the resulting graph has packing chromatic number $\pch(G) + 1$ or $c$ leaves were attached to $w$. In the latter case, continue by adding at most $c$ leaves to a new vertex. Proceeding in this way  we find a graph that has packing chromatic number $\pch(G) + 1$.
\qed

From Lemma~\ref{lem:abd}, we can now approach the question of determining if $(a,b,c)$ is realizable from a slightly different angle. Given positive integers $a$ and $b$, we define $m(a,b)$ to be the smallest integer such that $(a, b, m(a,b))$ is realizable. (Hence, $(a,b,c)$ is realizable if and only if $c\ge m(a,b)$.) We have already observed that $m(2,2) =2$ and it is easy to see that $m(a,a) =a$ for any $a\ge 2$. Indeed, this follows from the values of the invariants in complete graphs $K_a$, i.e., $\omega(K_a)=\chi(K_a)=\pch(K_a)=a$.

Next, we would like to study the relationship between $\chi(G)$ and $\pch(G)$ given an arbitrary graph $G$. As shown above, for any $b \ge 2$, we can find a graph where $\chi(G) = \pch(G) = b$. Is it possible that $(a, b, b)$ is realizable if $a<b$? This question was first considered in the seminal paper \cite{goddard-2008} where the following was shown.
\begin{proposition} {\rm (\cite[Proposition 2.6]{goddard-2008})}
For every graph $G$, if $\pch(G) = \chi(G)$, then $\omega(G) \ge \chi(G) - 2$.
\end{proposition}

Thus, if $(a, b, b)$ is realizable, then $a \ge b-2$. We further improve this, by showing that realizability of $(a,b,b)$ implies that $a=b$.

In the following proofs, we will be using the concept of {\em chromatic number criticality}. Recall that a graph $G$ is {\em $k$-critical} if $\chi(G)=k$ and for any proper subgraph $H$ of $G$, $\chi(H)<k$. It is well known that $k$-critical graphs are $k$-edge connected, and so the minimum degree $\delta(G)$ is at least $k-1$, cf.~\cite{west}.

\begin{theorem}\label{thm:chi=omega}
If $\chi(G) = \pch(G)\ge 3$, then $\omega(G) = \chi(G)$.
\end{theorem}

\proof
Let $k = \chi(G) = \pch(G)$. If $k=3$, then by \cite[Proposition 3.2]{goddard-2008}, $G$ contains the join of $K_2$ and an independent set as a subgraph. As the latter graph contains triangles, $\omega(G)=3$.

Now, let $G$ be a graph such that $\chi(G) = k=\pch(G)$, where $k\ge 4$. For the purpose of getting a contradiction, suppose that $\omega(G)<k$. We may assume that $G$ is a $k$-critical graph with respect to chromatic number. Indeed, if $G$ is not $k$-critical then it contains a proper subgraph $G'$, which is a $k$-critical graph. In particular,  $\chi(G')=k$, which in turn implies $\pch(G')=k$. Since $\omega(G')\le \omega(G)<k$, the non-existence of such a ($k$-critical) graph $G'$ would imply that also $G$ does not exist. Hence, we may assume that already $G$ is $k$-critical, and so $\delta(G) \ge k-1$.  It suffices to show the result is true for any connected graph $G$, hence we may, in addition, assume that $G$ is connected.

Let $c:V(G) \to [k]$ be a packing coloring of $G$ with color classes $V_1, \dots, V_k$. Since $V_i$ is an $i$-packing for each $i\in [k]$, the set $V_i$ is independent. This means that $(V_1, \dots, V_k)$ is a proper coloring with $k$($=\chi(G)$) colors. Therefore, there exists a vertex in each color class that is adjacent to a vertex of every other color. Furthermore, since every $x \in V_1$ has degree at most $k-1$ (otherwise $x$ would be adjacent to two vertices from some $V_i$, $i\ge 2$, which would then be at distance $2$) and yet $\delta(G) \ge k-1$, we know that $x$ is adjacent to exactly one vertex of colors $2, \dots, k$. Let $v_k \in V_k$ be a vertex of color $k$ that has a neighbor in every other color class. We let $v_i$, for each $3\le i \le k-1$, be the neighbor of $v_k$ with color $i$.

\medskip\noindent
{\underline{{\bf Claim.}} Vertex $v_k$ is the only vertex of $G$ with color $k$.}

\medskip\noindent
{\textbf{Proof.}} To see this, suppose that there exists another vertex $y \in V_k$ of color $k$. Since $G$ is connected, there exists a shortest $v_k,y$-path $P$ in $G$ of length at least $k+1$. We select $P = v_k w_1 w_2 w_3 w_4\cdots y$ such that $c(w_1)$ is smallest possible among all shortest $v_k,y$-paths.

Suppose first that $c(w_1) = 1$. As mentioned above, $w_1$ is adjacent to exactly one vertex of color $i$ for each $i$, $2 \le i \le k$. Thus, $w_1$ is adjacent to each $v_i$ for $3 \le i \le k-1$ as $d(w_1, v_i) \le 2$ for each $3 \le i \le k-1$. Indeed, otherwise a neighbor $x\ne v_i$ of $w_1$ of color $i$, $3 \le i \le k-1$, would be at distance at most $3$ from $v_i$.  It follows that $c(w_2) = 2$ and $c(w_3) = 1$ since $d(w_3, v_i) =3$ for each $3 \le i \le k$. This implies that $w_3$ is adjacent to $v_k$ since $k \ge 4$, which contradicts our choice of $P$. Thus, $c(w_1) >1$.

Next, assume that $c(w_1) =2$. Thus, $c(w_2) = 1$ as $d(w_2, v_i) \le 3$ for each $3 \le i \le k$, which also contradicts our choice of $P$ as $w_2$ is adjacent to $v_k$. Therefore, $w_1 = v_{\ell}$, $\ell >2$, and we know that $c(w_2) \in \{1, 2\}$ as $d(w_2, v_i) \le 3$ for each $3 \le i \le k$. If $c(w_2) = 1$, then $w_2$ is adjacent to $v_k$ as $d(v_k, w_2) = 2$. However, this contradicts our choice of $P$. Thus, we may assume $c(w_2) = 2$. Since $\delta(G)\ge k-1$, every vertex of color $s$, where $s>1$, has a neighbor of color 1.
In particular, $w_2$ has a neighbor of color $1$, call it $x$. It follows that $x$ is adjacent to $v_k$, meaning that $P' = v_k x w_2 w_3 w_4 \cdots y$ is a shortest path where $c(x) < c(w_1)$, contradicting our choice of $P$. Therefore, we may conclude that $v_k$ is the only vertex of color $k$. \smallqed

\medskip
Next, we claim that for each $\ell$, $3 \le \ell \le k-1$, $v_{\ell}$ is the only vertex of $G$ of color $\ell$. Indeed, fix $\ell$ and suppose there exists a vertex $y$ different from $v_{\ell}$ of color $\ell$. Since $G$ is connected, there exists a shortest $v_{\ell},y$-path $P$ in $G$. Among all such paths we select $P=v_{\ell} w_1w_2w_3 \cdots y$ such that $c(w_1)$ is as small as possible. Note that for each $i$, where $3 \le i \le k$, $d(w_1, v_i) \le 3$, since $v_{\ell}$ is adjacent to $v_k$ and $v_k$ is adjacent to $v_i$. Thus, $c(w_1) \in \{1, 2\}$. If $c(w_1) =1$, then $w_1$ is adjacent to $v_k$, meaning that $d(v_i, w_1) \le 2$ so $w_1$ is adjacent to each $v_i$ for $3 \le i \le k$. Therefore, $c(w_2) = 2$ and since $d(w_3, v_i) \le 3$, for each $3 \le i \le k$, $c(w_3) = 1$. Notice that $d(w_3, v_k) \le 4$ and since $k\ge 4$, $w_3$ is adjacent to $v_k$. However, this contradicts our choice of $P$ as $P' = v_{\ell}v_kw_3\cdots y$ is a shorter $v_{\ell},y$-path. Thus, $v_{\ell}$ is the only vertex of color $\ell$ for each $\ell \in \{3, \dots, k-1\}$.

Finally, we know that the graph induced by $\{v_3, \dots, v_k\}$ is a clique in $G$ since there exists a vertex of each color that is adjacent to all other colors. Furthermore, there exists a vertex $v_2 \in V_2$ that is adjacent to a vertex of every other color class. Let $v_1$ be a vertex of color $1$ that is adjacent to $v_2$. Thus, the graph induced by $\{v_2, \dots, v_k\}$ is a clique and since $d(v_1, v_i) \le 2$ for each $3 \le i \le k$, $v_1$ is adjacent to $v_i$ for each $2 \le i \le k$. It follows that $G$ contains a clique of size $k$, which is the final contradiction, implying that graphs $G$ with $\omega(G)<\chi(G)=\pch(G)=k$ do not exist.
\qed

By Theorem~\ref{thm:chi=omega}, $(2, 3, 3)$ is not realizable, and so $m(2, 3) > 3$. In fact, Theorem~\ref{thm:chi=omega} says that $(2, k, k)$ is not realizable for any $k \ge 3$ so we would like to compute $m(2, b)$ for any $b \ge 3$. An example of a graph $G$ that realizes $(2, 3, 4)$ is $C_5$, which is the Mycielskian of $K_2$. However, computing $m(2, b)$ becomes difficult rather quickly as $b$ gets larger. What we can say is that $m(2, b) \ge b+2$ for $b \ge 4$, as shown below.

\begin{theorem}
\label{thm:2-k-k+1}
If $k\ge 4$, then $(2, k, k+1)$ is not realizable.
\end{theorem}

\proof
Suppose there exists a graph $G$ of the form $(2, k, k+1)$ for some $k \ge 4$. Let $c:V(G) \to [k+1]$ be a $(k+1)$-packing coloring with color classes $V_1, \dots, V_{k+1}$. Let $H$ be the graph induced by $V_{k-2}\cup V_{k-1} \cup V_{k} \cup V_{k+1}$  and suppose $H$ is bipartite. This means we can properly color the vertices of $H$ using only two colors and in turn implies that $\chi(G) < k$, which is a contradiction. Therefore, $H$ is a triangle-free graph that contains odd cycles. Let $C = x_1x_2\cdots x_{n}$ be an odd cycle in $H$ of shortest length and note that in $G$, each vertex of $H$ is colored $k-2, k-1, k,$ or $k+1$. Thus, $C$ must be $(2, 3, 4, 5)$-packing colorable (i.e., $S$-packing colorable for the sequence $S=(2,3,4,5)$), which is not possible if $C\cong C_5$. Hence $C\cong C_n$, where $n\ge 7$ is an odd integer. As $C$ is $(2, 3, 4, 5)$-packing colorable, at most $\lfloor n/(i+1) \rfloor$ vertices can be assigned the color $i$ for each $2 \le i \le 5$. Let $(W_2,W_3,W_4,W_5)$ be a $(2, 3, 4, 5)$-packing coloring of $C$. Now
\[n=\sum_{i=2}^5 |W_i| \le \sum_{i=2}^5 \lfloor n/(i+1)\rfloor \le \frac{57}{60}n,\]
which is a contradiction. Hence, $C$ is not $(2, 3, 4, 5)$-packing colorable, which also implies that it is not $(k-2,k-1,k,k+1)$-packing colorable, for any $k\ge 4$, and thus no such graph $G$ exists.
\qed

Next, we improve Theorem~\ref{thm:2-k-k+1} by proving that $(2,k,k+2)$ is not realizable for any $k\ge 4$. We start with the case $k=4$.
\begin{lemma}
Triple $(2,4,6)$ is not realizable.
\label{lem:246}
\end{lemma}
\proof Suppose that there exists a graph $G$ such that $\omega(G)=2,\chi(G)=4$ and $\pch(G)=6$. Clearly, $G$ contains as a subgraph a 4-critical graph (with respect to chromatic number), say $H$, and we claim that $\pch(H)=6$. Indeed, since $H$ is triangle-free, $\omega(H)=2$, and by Theorem~\ref{thm:chi=omega}, $\pch(H)>\chi(H)=4$. By Theorem~\ref{thm:2-k-k+1}, $H$ cannot be a $(2,4,5)$ graph, hence $\pch(H)=6$. This implies that under the assumption that a $(2,4,6)$ graph exists, there are $(2,4,6)$ graphs that are 4-critical. This in turn implies that there exists a $(2,4,6)$ graph, say $G$, with $\delta(G)\ge 3$.

Consider a packing coloring $c$ of $G$ inducing a partition $(V_1,\ldots,V_6)$ into the color classes, where $V_i$ consists of the vertices that are assigned color $i$. As $\chi(G)=4$, the graph $G\setminus V_1$ cannot be bipartite, therefore it contains an odd cycle $C$. Clearly, $|V(C)|\ge 5$.

First, assume that there exist adjacent vertices $u,v\in V(C)$ such that $c(u)=2$ and $c(v)=3$. Let $x\ne u$ be the neighbor of $v$ in $C$, $y\ne v$ the neighbor of $u$ in $C$, and $z$ the other neighbor of $y$ in $C$. It is easy to see that $\{c(x),c(y),c(z)\}=\{4,5,6\}$. Now, as $\delta(G) \ge 3$, vertex $u$ has another neighbor in $G$, let it be $w$. By the distribution of the colors in $C$, we find that $c(w)$=1. Note that $w$ has at least two other neighbors in $G$, and they cannot be $v$ or $y$, since $G$ has no triangles. If one of these two neighbors is also different from $x$ and $z$, then we get a contradiction, because this vertex is then at distance at most $i$ from a vertex of color $i$ for every $i\in [6]$. Hence the only remaining possibility is that $w$ is adjacent to both $x$ and $z$. This implies that $|V(C)|> 5$ because $G$ has no triangles. Consider the neighbor $x'$ of $x$ in $C$, and its neighbor $x''$ in $C$ that is not $x$. Note that the only possibility for the color of $x'$ is that $c(x')=2$. In addition, the only possibility for the color of $x''$ is that $c(x'')=4$, which is possible only in the case when also $c(y)=4$. Now, consider the neighbor $z'$ of $z$, distinct from $y$ in $C$, and the neighbor $z''$ of $z'$ in $C$ that is distinct from $z$. There are two possibilities. Either $c(z')=3$ and $c(z'')=2$ , or $c(z')=2$ and $c(z'')=3$. In the first case note that vertices $z''$ and $z'$ are in an analogous setting as vertices $u$ and $v$. Since we deduced above that the neighbor of the neighbor of $v$ on $C$, which is vertex $x'$, must receive color $2$, we also infer that the neighbor of the neighbor of $z'$ on $C$, which is vertex $y$, must receive color $2$. This is a contradiction, because we already established that $c(y)=4$. Finally, if $c(z')=2$ and $c(z'')=3$, then as in the case of $u$ and $v$, we infer that $z'$ must have a neighbor $w'$ such that $c(w')=1$ and $w$ is adjacent to $y$ and also to the neighbor of $z''$ in $C$ different from $z'$; let it be called $t$. Now, $t$
is at distance at most $i$ from a vertex of color $i$ for every $i\in [6]$, which is the final contradiction, implying that no two adjacent vertices in $C$ can receive colors $2$ and $3$.

So the second case is that no two neighboring vertices in $C$ can receive colors $2$ and $3$. Therefore, as we pass along $C$, vertices with colors from $\{4,5,6\}$ appear one right after the other with possible gaps of at most one vertex (with color $2$ or $3$) in between. The longest possible pattern that vertices with colors from $\{4,5,6\}$ can form while passing along $C$ is $4-5-6-4-5$, where we did not write the vertices with color $2$ and $3$ that lie between them (we cannot continue the pattern with color $6$, because the next vertex is at distance $6$ from the vertex with color $6$ in the pattern). This implies that $|V(C)|\le 9$, hence $C$ can only be of length $5,7$ or $9$. It is easy to see that $C$ cannot have $9$ or $7$ vertices, since one can use only one vertex of each of the colors from $\{4,5,6\}$ (because the diameters of these two cycles are at most $4$) and we get in a contradiction with how the colors $2$ and $3$ are distributed in $C$.

Finally, we are left with the case that $|V(C)|=5$ and vertices with colors $2$ and $3$ are not adjacent in $C$. Clearly, the remaining vertices in $C$ receive colors $4,5,$ and $6$. Now, as $\delta(G)\ge 3$, the vertex $x$ with color $2$ has a neighbor $s$ in $G$, which is not in $C$, and it is obvious that $c(s)=1$. Vertex $s$ has at least two neighbors in $G$ besides $x$. Not both of these neighbors of vertex $s$ can be in $C$ because $G$ is triangle-free. A neighbor of $s$ that is not in $C$ can only receive color $3$, hence there can be only one such neighbor. The other neighbor of vertex $s$ must thus lie in $C$ and is not adjacent to a neighbor of $x$ on $C$. Noting that two vertices with color $3$ are at distance at most 3 we derive the final contradiction, by which the proof is complete.
\qed

\begin{theorem}
\label{thm:2-k-k+2}
Triple $(2,k,k+2)$ is not realizable for any $k$, $k\ge 4$.
\end{theorem}
\proof
The case $k=4$ was proven in Lemma~\ref{lem:246}, hence we may assume that $k\ge 5$. Suppose that there exists a graph $G$ such that $\omega(G)=2,\chi(G)=k$ and $\pch(G)=k+2$. Clearly, $G$ contains as a subgraph a $k$-critical graph (with respect to chromatic number), say $H$, and we claim that $\pch(H)=k+2$. Indeed, since $H$ is triangle-free, $\omega(H)=2$, and by Theorem~\ref{thm:chi=omega}, $\pch(H)>\chi(H)=k$. By Theorem~\ref{thm:2-k-k+1}, $H$ cannot be a $(2,k,k+1)$ graph, hence $\pch(H)=k+2$. This implies that under the assumption that a $(2,k,k+2)$ graph exists, there are $(2,k,k+2)$ graphs that are $k$-critical. This in turn implies that there exists a $(2,k,k+2)$ graph, say $G$, with $\delta(G)\ge k-1$.

Consider a packing coloring $c$ of $G$ inducing a partition $(V_1,\ldots,V_{k+2})$ into the color classes, where $V_i$ consists of the vertices that are assigned color $i$. As $\chi(G)=k$, the graph $G\setminus (V_1\cup\cdots\cup V_{k-3})$ cannot be bipartite, therefore it contains an odd cycle $C$. Clearly, $|V(C)|\ge 5$. Suppose $n=|V(C)|>5$. By the pigeon-hole principle, using also that the vertices with color $i$ must be more than distance $i$ apart, we infer that the number of vertices in $C$ having color $i$ is at most $\max\{1,\lfloor \frac{n}{i+1}\rfloor\}$. Altogether the number of vertices in $C$, by taking into account the available colors, is at most
$$\sum_{i=k-2}^{k+2}{\max\{1,\lfloor \frac{n}{i+1}\rfloor\}}.$$
As it turns out, this sum is strictly less than $n$, when $n\ge 7$. In particular, in the smallest case, where $k=5$ and $n=7$, we have
$\sum_{i=3}^{7}{\max\{1,\lfloor \frac{7}{i+1}\rfloor\}}=5.$

Finally, suppose that all odd cycles $C$ in $G\setminus (V_1\cup\cdots\cup V_{k-3})$ have length $5$. We restrict to the case when $k=5$. (As will be clear from the proof, the proof when $k>5$ follows similar lines only that a contradiction may appear even earlier.) Hence, since $k=5$, we are considering $G\setminus (V_1\cup V_2)$, and so the vertices in $C$ must get all colors from $3$ to $7$, each vertex a distinct color. Now, consider $G\setminus V_1$ and note that it need not be connected. If a component of $G\setminus V_1$ is isomorphic to $C_5$ in which vertices receive colors from $3$ to $7$, then this component is clearly 3-colorable. The same conclusion holds for a component of $G\setminus V_1$ whose subgraph induced by the vertices with colors from $3$ to $7$ is bipartite;  such a component of $G\setminus V_1$ is also 3-colorable. In either case this implies that $G$ is 4-colorable, which is a contradiction. Hence there must be a component in $G\setminus V_1$ such that to a 5-cycle $C$, in which vertices receive colors from $3$ to $7$, a vertex $x$ with color $2$ is attached as a neighbor of some vertex in $C$. Consider now this subgraph in $G$, and recall that $\delta(G)\ge 4$. Hence, $x$ has three more neighbors in $G$, at least two of which are not in $C$, because $G$ is triangle-free. In any case, regardless of how many neighbors $x$ has in $C$, we infer that two neighbors $x_1,x_2$ of $x$ (which are not in $C$) receive color $1$. Clearly, $x_1$ can have at most two neighbors in $C$ because $G$ is triangle-free. But if $x_1$ has two neighbors in $C$, then a neighbor $x'$ of $x_1$, which is not in $C$, is at most 3 apart from the vertex in $C$ with color $3$; as vertices with all other available colors are also too close to $x'$, we get a contradiction. Thus, $x_1$ can have at most one neighbor in $C$.
Now, if $x_1$ has a neighbor in $C$, then only color $3$ is possible for the other (at least two) neighbors of $x_1$, which gives us a contradiction. If on the other hand, $x_1$ is not adjacent to a vertex of $C$, then it has at least three neighbors, for which only colors $3$ and $4$ are available, which is the final contradiction. (Note that if $k>5$, available colors are bigger while distances in the subgraph are the same, therefore the last subcase involving $x_1$ gives an immediate contradiction.)
\qed

Summarizing the results concerning the function $m$, we first note that
$m(2,3)=4$ can be extended to an arbitrary $k$, $k\ge 2$, as follows. First, $m(k,k+1)>k+1$ by Theorem~\ref{thm:chi=omega}. On the other hand, $\omega(M(K_k))=k$, $\chi(M(K_k))=k+1$, and $\pch(M(K_k))=k+2$ (by Theorem~\ref{thm:mycielski}), which implies that $(k,k+1,k+2)$ is realizable for any $k$, $k\ge 2$. Combining both observations, we get $m(k,k+1)=k+2$.

\begin{table}[!ht]
\begin{tabular}{c|ccccccccc}
$a\backslash b$ &  2 & 3 & 4 & 5 & 6 & 7 & 8 & 9&10\\\hline
2 & 2 & 4 & 7 & & & & & &\\
3 & - & 3 & 5 & 6/9 &&&& &\\
4 & - & - & 4 & 6 & 7/11 &&&&\\
5 & - & - &- & 5& 7& 8/13 &&&\\
6 & - & - &- & - & 6& 8& 9/15&& \\
7 & - & - & - &- & -& 7& 9& 10/17&\\
8 & - & - & - & - &-& -& 8 & 10& 11/19 \\
\end{tabular}
\caption{The entry in row $a$ and column $b$ presents the known value of $m(a,b)$, while the entry separated by '$/$' present currently known lower and upper bound on $m(a,b)$.}
\label{tab:m}
\end{table}

Table~\ref{tab:m} summarizes the results on $m(a,b)$ presented so far. It can be complemented by the upper bound on the packing chromatic number of a graph obtained from the complete graphs by applying  the Mycielskian operation several times. Let us inductively define $M^k(G)$ as $M(M^{k-1}(G))$, where $M^1(G)$ is just the Mycielskian $M(G)$ of a graph $G$.
Applying Corollary~\ref{cor:mycielski-diam-2} inductively, starting from a complete graph, and using the fact that $\alpha(M(G))\ge |V(G)|$ for any graph $G$, we can prove by induction that
$$\pch(M^k(K_n))\le 2^{k-1}(n+1)+1,$$
for any $k\ge 1$. This implies that $m(n,n+k)\le 2^{k-1}(n+1)+1$ for any $k\ge 1$. In particular, $m(n,n+2)\le 2n+3$, which is used in Table~\ref{tab:m} as the upper bound values.

We suspect that the lower bound values in Table~\ref{tab:m} could be improved. After a close examination of the smallest case concerning the $(3,5,6)$ realizability, we pose the following conjecture.
\begin{conjecture}
There exists no graph $G$ with $\omega(G)=3, \chi(G)=5$ and $\pch(G)=6$. In other words, $(3,5,6)$ is not realizable.
\end{conjecture}
In fact, we suspect that $(k,k+2,k+3)$ might not be realizable for any $k\ge 3$.

\section{A lower bound on the packing chromatic number}
\label{sec:g-results}

\begin{proposition}
\label{prop:lowerbound}
If $G$ is a graph of maximum degree $\Delta(G)$, then
$$\pch(G)\ge \Delta(G)-\alpha(G)+2\,.$$
Equality is achieved if $\Delta(G)=n(G)-1$.
\end{proposition}

\proof
Let $r=\pch(G)$ and suppose that $(V_1,\ldots,V_r)$ is an $r$-packing coloring of $G$.  Let $x$ be any vertex in $G$.
If $x\in V_1$, then $x$ is adjacent to at most one vertex in $V_j$ for each $j\in [r] - \{1\}$.  Indeed if $x$ had two neighbors, say
$y_1$ and $y_2$, that both belong to $V_j$ for some $j \ge 2$, then $d_G(y_1,y_2)=2$, which contradicts the fact that $V_j$ is a $j$-packing.  This
implies that $\deg(x) \le r-1$.  If $x\in V_i$ for some $i \ge 2$, then $x$ has at most $|V_1|$ neighbors in $V_1$ and at most one neighbor in
$V_j$ for each $j\in[r]-\{1,i\}$.  Hence, $\deg_G(x) \le |V_1|+r-2 \le \alpha(G)+r-2$.  In both cases $\deg_G(x) \le \alpha(G)+r-2$, and it follows
that $\pch(G)\ge \Delta-\alpha(G)+2$.

Suppose that $G$ is a graph that has a vertex of degree $n(G)-1$. If $G$ is complete, the result is clear. Otherwise such a graph has diameter $2$ and thus from Proposition~\ref{prp:indep} we get $\pch(G)=n(G)-\alpha(G)+1=\Delta-\alpha(G)+2$.
\qed

Let ${\cal H}$ be the class of graphs $H$ constructed in the following way. Let $r\ge 3$ and $s\ge 2$ be positive integers. Let $A$ be a complete graph of order $r$ with three specified vertices $a_1,a_2$ and $a$.
Let $B$ be a complete graph of order $s$ with two specified vertices $b$ and $b_1$. Then let $H$ be any graph of order $r+s$ constructed from the disjoint union of $A$ and $B$ together
with a new vertex $z$ as follows.  Add an edge between $z$ and every vertex of $B\setminus \{b\}$ and then identify the vertices $a$ and $b$, call this vertex
$w$.  Any missing edge,
other than $a_1b_1$, can be added if it is not incident to $z$ or to $a_2$.

Note that in a graph $H\in {\cal H}$ as constructed above, the vertex $w$ is adjacent to every other vertex except $z$.
Let $V_1=\{a_1,b_1\}$, $V_2=\{z,a_2\}$ and for each $j$ such that $3\le j \le r+s-2$, let $V_j$ be a single vertex.  It is easy to verify that
$(V_1,V_2,\ldots,V_{r+s-2})$ is a packing coloring, and indeed that $\pch(H)=r+s-2$.

\begin{theorem}
If $G$ is a graph with $\alpha(G) = 2$, then $\chi_{\rho}(G) =  \Delta(G) - \alpha(G) + 2$ if and only if $\Delta(G) = n(G)-1$ or $G\in {\cal H}$.
\end{theorem}

\proof
Let $G$ be a graph such that $\alpha(G) = 2$.  If $\Delta(G) = n(G)-1$, then the diameter of $G$ is 2 and from Proposition~\ref{prp:indep} it follows that
$\chi_{\rho}(G) =  n(G) - \alpha(G) + 1=\Delta(G) - \alpha(G) + 2$.  If  $G\in {\cal H}$, then by the paragraph just before the theorem we see that
$\chi_{\rho}(G) =  \Delta(G) - \alpha(G) + 2$.

For the converse let $r = \chi_{\rho}(G)$ and suppose $c:V(G) \to [r]$ is a packing coloring with color classes $V_1, \dots, V_r$. Let $w$ be a vertex of $G$ with degree $\Delta(G)=\chi_{\rho}(G) + \alpha(G) - 2 = \chi_{\rho}(G)$. If $w \in V_1$, then $w$ has at most one neighbor in $V_i$ for each $i \ge 2$, which implies that $\deg(w) \le r-1$ and  contradicts the assumption that $\Delta(G) = r$. Thus, we may assume that $w \not\in V_1$. Since $\deg(w) = \chi_{\rho}(G)$, we know that $w$ is adjacent to exactly one vertex in $V_j$ for each $j \not\in \{1, c(w)\}$ and to two vertices in $V_1$. Moreover, $V_1$ is a maximum independent set, meaning that every vertex of $V(G) \setminus V_1$ is adjacent to some vertex in $V_1$. We write $V_1 = \{v_1, v_2\}$ and let $y_i \in V_i$ be the vertex adjacent to $w$ for each $i \not\in \{1, c(w)\}$. Note that if $V(G) = N[w]$, then $\Delta(G) = n(G)-1$ and we are done. So we shall assume that $N[w] \not= V(G)$.

Let $z \in V(G) \setminus N[w]$ and observe that $z \not\in V_{c(w)}$ and $z \not\in V_1$. Thus, $z \in V_{c(z)}$ where  $c(z) \not\in \{1, c(w)\}$. As stated above, $z$ is adjacent to some $v_i$ in $V_1$. Without loss of generality, we may assume that $z$ is adjacent to $v_1$. It follows that $zv_1wy_{c(z)}$ is a path of length $3$ in $G$ which cannot exist if $c(z) > 2$. Hence, $c(z) = 2$ and  $V_2 = \{z, y_2\}$ since $\alpha(G)=2$.  Finally, we point out that $|V_i| = 1$ for all $i \not\in \{1, 2\}$ and $\chi_{\rho}(G) = n(G)-2$.

Note that $z$ and $y_2$ are not adjacent but every other vertex of $G$ is adjacent to exactly one of $z$ or $y_2$.  Therefore, the two sets $N[y_2]$ and
$N[z]$ partition $V(G)$.  We claim that $N[z]$ and $N[y_2]$ are both complete subgraphs.  Let $u$ and $v$ be any two vertices of $N[y_2]$.  The vertex
$z$ is adjacent to neither $u$ nor $v$, and since $\alpha(G)=2$, it follows that $uv \in E(G)$ and hence $N[y_2]$ is complete.  Similarly, $N[z]$ is
a complete subgraph.  Referring to the description in the paragraph before the statement of the theorem, we can now complete the proof by letting
$A=N[y_2]$ with specified vertices $a_1=v_2$, $a_2=y_2$ and $a=w$.  Let $B=N(z)$ with specified vertices $b_1=v_1$ and $b=w$.  Hence $G\in {\cal H}$.
\qed

We conclude the paper by noting that a family of graphs $G$ from $\cal H$, which are obtained from the complete graph $K_n$ by removing the edges of a subgraph isomorphic to $K_{1,r}$, where $r+1<n$, has the property that $\omega(G)=\chi(G)=\pch(G)=\Delta(G)=n-1$. Thus this is another infinite family of graphs that realizes the triple $(n-1,n-1,n-1)$ for any $n$, for $n\ge 3$.

\section*{Acknowledgements}

This work was supported in part by the Ministry of Science of Slovenia
under the grant ARRS-BI-US/16-17-013.
B.B.\ and S.K.\ acknowledge the financial support from the Slovenian Research Agency (research core funding No.\ P1-0297) and that the project (Combinatorial Problems with an Emphasis on Games, N1-0043) was financially supported by the Slovenian Research Agency.
D.F.R.\ is supported by a grant from the Simons Foundation (Grant Number \#209654 to Douglas F. Rall).


\end{document}